\documentclass[12pt]{article}
\usepackage{graphicx}
\usepackage{amsmath,amsthm,amssymb,enumerate}
\usepackage{euscript,mathrsfs}
\usepackage{color}
\usepackage{dsfont}
\usepackage[left=2cm,right=2cm,top=3.5cm,bottom=3.5cm]{geometry}
\usepackage{color}
\usepackage[framemethod=tikz]{mdframed}
\allowdisplaybreaks

\usepackage{soul}
\usepackage{esint}
\catcode`\@=11 \@addtoreset{equation}{section}

\catcode`\@=12

\allowdisplaybreaks

\newtheorem{Theorem}{Theorem}[section]
\newtheorem{Proposition}[Theorem]{Proposition}
\newtheorem{Lemma}[Theorem]{Lemma}
\newtheorem{Corollary}[Theorem]{Corollary}

\theoremstyle{definition}
\newtheorem{Definition}[Theorem]{Definition}

\newtheorem{Remark}[Theorem]{Remark}

\newcommand{\bTheorem}[1]{
\begin{Theorem} \label{T#1} }
\newcommand{\eT}{\end{Theorem}}

\newcommand{\bProposition}[1]{
\begin{Proposition} \label{P#1}}
\newcommand{\eP}{\end{Proposition}}

\newcommand{\bLemma}[1]{
\begin{Lemma} \label{L#1} }
\newcommand{\eL}{\end{Lemma}}

\newcommand{\bCorollary}[1]{
\begin{Corollary} \label{C#1} }
\newcommand{\eC}{\end{Corollary}}

\newcommand{\bRemark}[1]{
\begin{Remark} \label{R#1} }
\newcommand{\eR}{\end{Remark}}

\newcommand{\bDefinition}[1]{
\begin{Definition} \label{D#1} }
\newcommand{\eD}{\end{Definition}}

\newcommand{\bfphi}{\boldsymbol{\varphi}}

\newcommand{\bFormula}[1]{
\begin{equation} \label{#1}}
\newcommand{\eF}{\end{equation}}

\newcommand{\Ov}[1]{\overline{#1}}

\newcommand{\TD}{\mathbb{T}^2}

\newcommand{\vr}{\varrho}

\newcommand{\tvr}{\tilde \vr}
\newcommand{\tvu}{{\tilde \vu}}
\newcommand{\tvt}{\tilde \vt}

\newcommand{\vt}{\vartheta}
\newcommand{\vu}{\vc{u}}

\newcommand{\vc}[1]{{\bf #1}}

\newcommand{\Div}{{\rm div}_x}
\newcommand{\Grad}{\nabla_x}

\newcommand{\dx}{\,{\rm d} {x}}

\newcommand{\dt}{\,{\rm d} t }

\newcommand{\intTD}[1]{\int_{\mathbb{T}^2} #1 \ \dx}

\newcommand{\intTS}[1]{\int_{\mathbb{R} \times \mathbb{T}^1} #1 \ \dx}

\newcommand{\ep}{\varepsilon}

\newcommand{\br}{ \nonumber \\ }

\def\softd{{\leavevmode\setbox1=\hbox{d}%
          \hbox to 1.05\wd1{d\kern-0.4ex{\char039}\hss}}}
\definecolor{Cgrey}{rgb}{0.85,0.85,0.85}
\definecolor{Cblue}{rgb}{0.50,0.85,0.85}
\definecolor{Cred}{rgb}{1,0,0}
\definecolor{fancy}{rgb}{0.10,0.85,0.10}

\newcommand\Cbox[2]{%
    \newbox\contentbox%
    \newbox\bkgdbox%
    \setbox\contentbox\hbox to \hsize{%
        \vtop{
            \kern\columnsep
            \hbox to \hsize{%
                \kern\columnsep%
                \advance\hsize by -2\columnsep%
                \setlength{\textwidth}{\hsize}%
                \vbox{
                    \parskip=\baselineskip
                    \parindent=0bp
                    #2
                }%
                \kern\columnsep%
            }%
            \kern\columnsep%
        }%
    }%
    \setbox\bkgdbox\vbox{
        \color{#1}
        \hrule width  \wd\contentbox %
               height \ht\contentbox %
               depth  \dp\contentbox
        \color{black}
    }%
    \wd\bkgdbox=0bp%
    \vbox{\hbox to \hsize{\box\bkgdbox\box\contentbox}}%
    \vskip\baselineskip%
}

\mdfdefinestyle{MyFrame}{%
	linecolor=black,
	outerlinewidth=1pt,
	roundcorner=5pt,
	innertopmargin=\baselineskip,
	innerbottommargin=\baselineskip,
	innerrightmargin=10pt,
	innerleftmargin=10pt,
	backgroundcolor=white!20!white}

\date{}

\begin{document}


\title{Glimm's method and density of wild data for the Euler system of gas dynamics}

\author{Elisabetta Chiodaroli 
\and Eduard Feireisl \thanks{The work of E.F. was partially supported by the
Czech Sciences Foundation (GA\v CR), Grant Agreement
21--02411S. The Institute of Mathematics of the Academy of Sciences of
the Czech Republic is supported by RVO:67985840. } 
}

\date{}

\maketitle

\medskip

\centerline{Dipartimento di Matematica, Universit\` a di Pisa,} 
	
\centerline{ Via F. Buonarroti 1/c, 56127, Pisa, Italy}

\medskip

\centerline{Institute of Mathematics of the Academy of Sciences of the Czech Republic}

\centerline{\v Zitn\' a 25, CZ-115 67 Praha 1, Czech Republic}

\date{}

\maketitle

\begin{abstract}
	
	We adapt Glimm's approximation method to the framework of convex integration to show density of wild data  
	for the (complete) Euler system of gas dynamics. The desired infinite family of entropy admissible solutions 
	emanating from the same initial data is 
	obtained via convex integration 
	of suitable Riemann problems pasted with local smooth solutions. In addition, the wild data belong to BV class.

\end{abstract}

\bigskip

{\bf Keywords:} Euler system of gas dynamics, wild data, Riemann problem, convex integration

\bigskip

\section{Introduction}
\label{i}

The celebrated Glimm method uses solutions of the Riemann problem as building blocks to construct weak solutions to non--linear conservation laws. 
More recently, the Riemann type initial data gave rise to several examples of ill--posedness of problems in fluid dynamics via the method of convex integration, see \cite{ChiDelKre}, 
\cite{ChiKre}, Markfelder and Klingenberg \cite{MarKli}, or, more recently, Al Baba et al \cite{ABKlKrMaMa}, Klingenberg et al. 
\cite{KlKrMaMa} to name only a few. The concept of \emph{wild data} emerged to identify the initial state of a system that leads to multiple solutions that are still physically admissible. 
Sz\' ekelyhidi and Wiedemann \cite{SzeWie} showed that the wild data generating infinitely many energy dissipating solutions of the \emph{incompressible} Euler system are dense in the $L^2$-topology.
A comparable result for the \emph{isentropic} Euler system was obtained by Chen, Vasseur, and Yu \cite{ChVaYu}. The admissibility criterion used in \cite{ChVaYu} based on the stipulation that the energy 
of the system never exceeds the initial energy may be rough and still compatible with non--physical (increasing) energy profiles. The problem of density of wild data for the isentropic Euler 
system was revisited in \cite{ChFe2022} in the class of physically admissible entropy solutions satisfying the energy inequality in the differential form on a possibly short time interval. Pursuing 
a similar strategy, we address the same problem in the context of the (complete) Euler system of gas dynamics.

\subsection{Euler system}

The Euler system of gas dynamics describes the evolution of the density $\vr= \vr(t,x)$, the temperature $\vt = \vt(t,x)$ and the velocity 
$\vu = \vu(t,x)$ by a system of field equations
\begin{align}
	\partial_t \vr + \Div (\vr \vu) &= 0, \label{E1}\\ 
	\partial_t (\vr \vu) + \Div (\vr \vu \otimes \vu) + \Grad p(\vr, \vt) &= 0, \label{E2}\\
	\partial_t \left( \frac{1}{2} \vr |\vu|^2 + \vr e(\vr, \vt) \right) + 
	\Div \left[ \left( \frac{1}{2} \vr |\vu|^2 + \vr e(\vr, \vt) + p(\vr, \vt) \right) \vu \right] &=0 \label{E3}
\end{align}
that express the physical principles of conservation of mass, linear momentum, and energy.
For definiteness, the pressure $p=p(\vr, \vt)$ and the internal energy $e=e(\vr, \vt)$ are chosen to obey the standard Boyle--Mariotte law,
\begin{equation} \label{E4} 
	p(\vr, \vt) = \vr \vt,\ e(\vr, \vt) = c_v \vt,\ c_v > 1.
	\end{equation}
For technical reasons, we consider the space periodic boundary conditions identifying the physical fluid domain with the flat torus 
\[
\Omega = \mathbb{T}^2 = \left\{ (x_1,x_2) \ \Big| \ x_1 \in [0,1]|_{\{0;1\}},\ x_2 \in [0,1]|_{\{0;1\}} \right\} .
\]
Possible generalizations are discussed in the concluding Section \ref{D}. The problem is completed by prescribing the initial conditions
\begin{equation} \label{E5}
\vr(0, \cdot) = \vr_0 ,\ \vt(0, \cdot) = \vt_0,\ \vu(0, \cdot) = \vu_0.
\end{equation}

As is well known, see e.g. the monograph by Dafermos \cite{D4a}, solutions of the Euler system may develop singularities in the form of shock waves even if the initial data are smooth. 
Weak (distributional) solutions have been introduced to capture the behaviour of discontinuities along with several \emph{admissibility criteria} to restore uniqueness in this 
larger class. In particular, the \emph{physically admissible solutions} to the Euler system are required to satisfy the entropy inequality
\begin{equation} \label{entro}
	\partial_t (\vr s(\vr, \vt)) + \Div( \vr s(\vr, \vt) \vu) \geq 0,\ \mbox{where}\ 
	s(\vr, \vt) = c_v \log \vt - \log \vr, 
\end{equation}
in the sense of distributions.

\subsection{``Wild'' solutions/data?} 

In the last decade, the method of convex integration revealed a number of rather unexpected facts concerning well/ill posedness of problems in fluid mechanics in the framework 
of weak solutions, see e.g. Buckmaster et al. \cite{BuDeSzVi}, \cite{BucVic} and the references cited therein. In particular, the complete Euler system 
\eqref{E1}--\eqref{E3}, with the initial data \eqref{E5} is ill--posed in the class of admissible weak solutions satisfying, in addition, the entropy inequality \eqref{entro}, 
see \cite{FeKlKrMa} for the periodic or impermeable bounded domain, and \cite{ABKlKrMaMa}, \cite{KlKrMaMa} for the Riemann problem. In particular, the Euler system 
admits \emph{wild data} - the initial conditions that give rise to infinitely many physically admissible solutions. 

Nonetheless, the fact that the energy \emph{equality} \eqref{E3} is included in the system may give some hope that wild data may be somehow ``exceptional'' for the 
complete Euler system, see \cite{FeKlMark}. In this paper, we show that the wild data for the Euler system with periodic boundary conditions are dense in the 
$L^p-$topology. In addition, we show certain regularity of the emerging solutions as well as density of the set on which non--uniqueness holds. In particular, 
\begin{itemize} 
	\item the wild data we construct are piecewise smooth, in particular in BV class;
	\item the corresponding solutions are differentiable with the exception of a set of small measure; 
	\item nonuniqueness  occurs in an ``almost'' dense subset of the physical domain.
	\end{itemize}

Our approach is based on a simple idea motivated by Glimm's approximation technique. First, we divide the physical space into a finite number of components, where we solve 
a Riemann type problem and make use of the ill--posedness results established in \cite{ABKlKrMaMa}, \cite{KlKrMaMa}. On the complementary part of the physical space, 
we consider the local in time strong solution. Thanks to the principle of finite speed of propagation, we paste the solution pieces together obtaining the desired result.

\section{Solvability of the Euler system}

We review some known facts concerning solvability of the Euler system in the class of both weak and strong solutions.

\subsection{Admissible (entropy) solutions}

\begin{Definition}[{\bf Admissible weak solution}] \label{DE1}
	A triple $(\vr, \vt, \vu)$ is called \emph{admissible weak solution} of the Euler system in $[0,T) \times \mathbb{T}^2$ emanating from the initial 
	data $(\vr_0, \vt_0, \vu_0)$ if the following holds:
	\begin{equation} \label{E6}
		\int_0^T \intTD{ \Big( \vr \partial_t \varphi + \vr \vu \cdot \Grad \varphi \Big) } \dt = - \intTD{ \vr_0 \varphi (0, \cdot) }
		\end{equation}
	for any $\varphi \in C^1_c([0,T) \times \TD)$; 
	\begin{equation} \label{E7}
		\int_0^T \intTD{ \Big( \vr \vu \cdot \partial_t \bfphi + \vr \vu \otimes \vu: \Grad \bfphi + p(\vr, \vt) \Div \bfphi \Big) } \dt = - \intTD{ \vr_0 \vu_0 \cdot \bfphi (0, \cdot) }	
		\end{equation}
	for any $\bfphi \in C^1_c([0,T) \times \TD; R^2)$;
	\begin{align} 
	\int_0^T &\intTD{ \left( \left( \frac{1}{2} \vr |\vu|^2 + \vr e(\vr, \vt) \right) \partial_t \varphi + 
	\left( \frac{1}{2} \vr |\vu|^2 + \vr e(\vr, \vt) + p(\vr, \vt) \right) \vu  \cdot \Grad \varphi \right) } \dt \br &\quad = - \intTD{
\left( \frac{1}{2} \vr_0 |\vu_0|^2 + \vr_0 e(\vr_0, \vt_0) \right) \varphi (0, \cdot) }
\label{E8}
	\end{align}
	for any $\varphi \in C^1_c([0,T) \times \TD)$; and	
\begin{equation} \label{E9}	
\int_0^T \intTD{ \Big( \vr s(\vr, \vt) \partial_t \varphi + \vr s(\vr, \vt) \vu \cdot \Grad \varphi \Big) } \dt \leq - 
\intTD{ \vr_0 s(\vr_0, \vt_0) \varphi (0, \cdot) }	
	\end{equation}
for any  $\varphi \in C^1_c([0,T) \times \TD)$, $\varphi \geq 0$.
	\end{Definition}

Despite the large number of recent results concerning existence of (infinitely many) weak solutions to the Euler 
system, the existence of an admissible solution for arbitrary, say bounded, initial data is still an open problem. 

\subsection{Riemann problem on a strip}
\label{R}

A special class of data that plays a crucial role in our analysis are the piece--wise constant so called Riemann 
data. Consider a spatial domain 
\[
\mathbb{R} \times \mathbb{T}^1 = \left\{ (x_1,x_2) \ \Big|\ x_1 \in \mathbb{R},\ x_2 \in [0,1]|_{\{0;1\}} \right\}
\]
an infinite 2-D strip.
The Riemann data are determined by constant vectors 
\begin{equation} \label{R1}
	(\vr_\ell, \vt_\ell, \vu_\ell),\ (\vr_r, \vt_r, \vu_r) \in (0,\infty) \times (0, \infty) \times \mathbb{R}^2 
\end{equation}	

\begin{Definition}[{\bf Riemann solution}] \label{DR1}
	A triple $(\vr_R, \vt_R, \vu_R)$ is called \emph{Riemann solution} of the Euler system in $[0,T) \times \mathbb{R} \times \mathbb{T}^1$ emanating from the Riemann data 
	data $(\vr_\ell, \vt_\ell, \vu_\ell)$, $(\vr_r, \vt_r, \vu_r)$ if the following holds:
	There exists $\lambda > 0$ such that 
	\begin{equation} \label{R1a}
	(\vr_R, \vt_R, \vu_R)(t,x) = (\vr_\ell, \vt_\ell, \vu_\ell) \ \mbox{if}\ x_1 < - \lambda t,\ 
	(\vr_R, \vt_R, \vu_R)(t,x) = (\vr_r, \vt_r, \vu_r) \ \mbox{if}\ x_1 > \lambda t;
	\end{equation}
	\begin{equation} \label{R2}
		\int_0^T \intTS{ \Big( \vr_R \partial_t \varphi + \vr_R \vu_R \cdot \Grad \varphi \Big) } \dt = - \intTS{ \left( \mathds{1}_{x_1 < 0} \vr_\ell + 
		\mathds{1}_{x_1 > 0} \vr_r \right)	 \varphi (0, \cdot) }
	\end{equation}
	for any $\varphi \in C^1_c([0,T) \times (\mathbb{R} \times \mathbb{T}^1))$; 
	\begin{align} \label{R3}
		\int_0^T & \intTS{ \Big( \vr \vu \cdot \partial_t \bfphi + \vr \vu \otimes \vu: \Grad \bfphi + p(\vr, \vt) \Div \bfphi \Big) } \dt \br 
		&= - \intTS{ \left( \mathds{1}_{x_1 < 0} \vr_\ell \vu_\ell + 
			\mathds{1}_{x_1 > 0} \vr_r \vu_r	    \right)\cdot \bfphi (0, \cdot) }	
	\end{align}
	for any $\bfphi \in C^1_c([0,T) \times (\mathbb{R} \times \mathbb{T}^1) ; \mathbb{R}^2)$;
	\begin{align} 
		\int_0^T &\intTS{ \left( \left( \frac{1}{2} \vr |\vu|^2 + \vr e(\vr, \vt) \right) \partial_t \varphi + 
			\left( \frac{1}{2} \vr |\vu|^2 + \vr e(\vr, \vt) + p(\vr, \vt) \right) \vu  \cdot \Grad \varphi \right) } \dt \br &= - \intTS{
			\left[ \mathds{1}_{x_1 < 0} \left( \frac{1}{2} \vr_\ell |\vu_\ell|^2 + \vr_\ell e(\vr_\ell, \vt_\ell) \right) + \mathds{1}_{x_1 > 0} \left( \frac{1}{2} \vr_r |\vu_r|^2 + \vr_r
			 e(\vr_r, \vt_r) \right) 
			 \right] \varphi (0, \cdot) }
		\label{R4}
	\end{align}
	for any $\varphi \in C^1_c([0,T) \times(\mathbb{R} \times \mathbb{T}^1) )$; 	
	\begin{align} 	
		\int_0^T &\intTS{ \Big( \vr s(\vr, \vt) \partial_t \varphi + \vr s(\vr, \vt) \vu \cdot \Grad \varphi \Big) } \dt \br &\leq - 
		\intTS{ \left( \mathds{1}_{x_1 < 0} \vr_\ell s(\vr_\ell, \vt_\ell) +   \mathds{1}_{x_1 > 0} \vr_r s(\vr_r, \vt_r)\right) \varphi (0, \cdot) }
		\label{R5}	
	\end{align}
	for any  $\varphi \in C^1_c([0,T) \times (\mathbb{R} \times \mathbb{T}^1))$, $\varphi \geq 0$.

\end{Definition}

The following result was proved by Klingenberg et al \cite[Theorem 1.1]{KlKrMaMa}, see also Al Baba et al. \cite{ABKlKrMaMa}. 

\begin{Proposition} \label{PR1} {\bf (Ill posedness for the Riemann problem).}
	
	\noindent	
	There exist Riemann data $(\vr_\ell, \vt_\ell, \vu_\ell)$, $(\vr_r, \vt_r, \vu_r)$ such that the Riemann problem for the Euler system 
	admits infinitely many solutions $(\vr_R, \vt_R, \vu_R)$ in $[0,T) \times (\mathbb{R} \times \mathbb{T}^1)$, $T > 0$ arbitrary. Moreover, all 
	solutions satisfy \eqref{R1a} with the same constant $\lambda$, and there are constants 
	\[
	0 < \underline{\vt} \leq \Ov{\vt},\ 
	0 < \underline{\vr} \leq \Ov{\vr},\ \Ov{\vu} > 0 
	\]
	such that 
	\begin{equation} \label{R7}
	\underline{\vt} \leq \vt_R \leq \Ov{\vt},\ 	\underline{\vr} \leq \vr_R \leq \Ov{\vr},\ |\vu_R| \leq \Ov{\vu}
		\end{equation}
a.e. in $(0,T) \times (\mathbb{R} \times \mathbb{T}^1)$.	
	\end{Proposition}

\begin{Remark} \label{RR1}
	As a matter of fact, the result in \cite{KlKrMaMa} is stated in the spatial domain $\mathbb{R}^2$. However, exactly the same 
	arguments yield the desired result on the strip $\mathbb{R} \times \mathbb{T}^1$. To see this, it is enough to observe that:
	\begin{itemize}
		\item the subsolutions constructed in \cite{KlKrMaMa} are piecewise constant and independent of $x_2$, in particular periodic in 
		$x_2$; 
		\item the oscillatory lemma (Lemma 2.6 in \cite{KlKrMaMa}) holds on arbitrary domain.
		\end{itemize} 
\end{Remark}

\subsection{Local existence of regular solutions}
\label{L}

Finally, we recall the well--known results concerning the existence of local--in--time 
strong solutions emanating from regular initial data, see e.g. Benzoni-Gavage and Serre 
\cite[Chapter 13, Theorem 13.1]{BenSer}. 

\begin{Proposition}[{\bf Local existence for regular data}] \label{PL1}
Let the initial data 
\[
	\vr_0 \in W^{k,2}(\TD),\ \inf_{\TD} \vr_0 > 0,\ \vt_0 \in W^{k,2}(\TD),\ \inf_{\TD} \vt_0 > 0,\ \vu_0 \in W^{k,2}(\TD; \mathbb{R}^2),\ k > 2,
\]
be given. 

Then there exists $T_{\rm max} > 0$ depending only on the norm of the data in the aforementioned spaces such that 
the Euler system admits a classical (continuously differentiable) solution $(\vr, \vt, \vu)$ in 
$[0, T_{\rm max}) \times \TD$, unique in the class 
\[
\vr, \vt \in C([0,T]; W^{k,2}(\TD)),\ \vu \in C([0,T]; W^{k,2}(\TD; \mathbb{R}^2)) 
\]
for any $0 < T < T_{\rm max}$.
	
	\end{Proposition}

\section{Main result}
\label{M}

\begin{Definition}[{\bf Wild data}] \label{DM1}
	The initial data $(\vr_0, \vt_0, \vu_0)$ are called \emph{wild} if there exists a time $T > 0$ such that the Euler system 
	admits infinitely many admissible weak solutions emanating from $(\vr_0, \vt_0, \vu_0)$ in $[0,\tau] \times \TD$ for any $0 < \tau < T$.

	\end{Definition}

\begin{Theorem}[{\bf Density of wild data}] \label{TM1}
	Let 
	\begin{equation} \label{M1}
	\vr_0 \in W^{k,2}(\TD),\ \inf_{\TD} \vr_0 > 0,\ \vt_0 \in W^{k,2}(\TD),\ \inf_{\TD} \vt_0 > 0,\ \vu_0 \in W^{k,2}(\TD; \mathbb{R}^2),\ k > 2,
		\end{equation} 
and $1 \leq q < \infty$ be given.

Then for any $\ep > 0$, there exist initial data $(\vr_{0, \ep}, \vt_{0, \ep}, \vu_{0, \ep})$ enjoying the following properties:
\begin{itemize}
\item
The data	
$(\vr_{0, \ep}, \vt_{0, \ep}, \vu_{0, \ep})$ are piecewise smooth, specifically, there are finitely many points $x_1^1, \dots, x^N_1$ such that
$(\vr_{0, \ep}, \vt_{0, \ep}, \vu_{0, \ep})(x_1,x_2)$ are continuously differentiable whenever $x_1 \not \in \{ x_1^1, \dots, x^N_1 \}$, 
$x_2 \in \mathbb{T}^1$.
\item 
There exists $T > 0$ such that the Euler system admits infinitely many admissible weak solutions $(\vr^n, \vt^n, \vu^n)_{n \in N}$ 
in $L^\infty([0,T) \times \TD; \mathbb{R}^4)$ emanating from the initial data $(\vr_{0, \ep}, \vt_{0, \ep}, \vu_{0, \ep})$ such that 
\begin{equation} \label{M2}
(\vr^n, \vt^n, \vu^n)|_{[0,\tau) \times B_\ep} \not \equiv (\vr^m, \vt^m, \vu^m)|_{[0,\tau) \times B_\ep}, \forall \ m \ne n	
	\end{equation}
whenever $0 < \tau < T$ and $B_\ep \subset \TD$ is a ball of radius $\ep$.
\item 
\begin{equation} \label{M3}
	\left\| \Big(\vr_{0,\ep} - \vr_0;  \vt_{0,\ep} - \vt_0; \vu_{0,\ep} - \vu_0 \Big) \right\|_{L^q(\TD; \mathbb{R}^4)} \leq \ep;
	\end{equation}

\item 
\begin{equation} \label{M4}
\inf_{[0,T) \times \TD} \vr^n > 0,\ \inf_{[0,T) \times \TD}	 \vt^n > 0, 
	\end{equation}
and there exists a compact set $\mathcal{S} \subset  \TD$, $|\mathcal{S}| < \ep$ such that 
$(\vr^n, \vt^n, \vu^n)$ are continuously differentiable in $(0,T) \times (\TD \setminus \mathcal{S})$ for any $n \in N$.
\end{itemize}

\end{Theorem}

\begin{Corollary} \label{CM1}
	
	Wild data are dense in $L^q(\TD; [0, \infty)^2 \times \mathbb{R}^2)$ for any finite $q \geq 1$.
	
	\end{Corollary}

\section{Proof of the main result}
\label{P}

\subsection{Glimm's partition}

Fix $N = N(\ep)$ distinct points $x^1_1 < \dots < x^N_1$, $x^i_1 \in \mathbb{T}^1$, $i=1,\dots,N$ such that 
\begin{equation} \label{P1}
|x^{i+1}_1 - x^i_1| < \ep \ \mbox{for all}\ i = 1,\dots, N, 
\end{equation}
where we have identified $x^{N+1}_1 = x^1_1$. Accordingly, for any (open) ball $B_\ep$ of radius $\ep$, there exists $m$ such that 
\begin{equation} \label{P2}
	B_\ep \cap \left\{ (x^m_1, x_2) \Big| \ x_2 \in \mathbb{T}^1 \right\} \ne \emptyset.
\end{equation}	

\subsection{Initial data, approximation} \label {Papprox}

Let $(\vr_0, \vt_0, \vu_0)$ be given as in \eqref{M1}. For $\delta > 0$ small enough, we find $(\vr_{0, \delta}, \vt_{0, \delta}, \vu_{0, \delta})$ 
such that 
\begin{equation} \label{P4}
\vr_{0,\delta}(x_1, x_2) = \vr_0(x_1,x_2),\ 
\vt_{0,\delta}(x_1, x_2) = \vt_0(x_1,x_2),\  \vu_{0,\delta}(x_1, x_2) = \vu_0(x_1,x_2)	
\end{equation}	
for any $x = (x_1,x_2)$ in the set
\begin{equation} \label{P5}
x_1,\ |x_1 - x^i_1| > 3 \delta.\ \mbox{for all}\ i = 1,\dots, N,\ x_2 \in \mathbb{T}^1,
\end{equation}
where the points $x^i_1$ have been introduced in \eqref{P1}.

Let $(\vr_\ell, \vt_\ell, \vu_\ell)$, $(\vr_r, \vt_r, \vu_r)$, be the Riemann data from Proposition \ref{PR1} giving rise to infinitely many admissible weak solutions 
to the Euler system. We set 
\begin{align} 
	\vr_{0,\delta} (x_1, x_2) &= \vr_\ell \ \mbox{if}\ x^i_1 - 2 \delta \leq x_1 \leq x^i_1 - \delta ,\ 
	\vr_{0,\delta} (x_1, x_2) = \vr_r \ \mbox{if}\ x^i_1 + \delta \leq x_1 \leq x^i_1 + 2 \delta,\ x_2 \in \mathbb{T}^1, \br
		\vt_{0,\delta} (x_1, x_2) &= \vt_\ell \ \mbox{if}\ x^i_1 - 2 \delta \leq x_1 \leq x^i_1 - \delta ,\ 
	\vt_{0,\delta} (x_1, x_2) = \vt_r \ \mbox{if}\ x^i_1 + \delta \leq x_1 \leq x^i_1 + 2 \delta,\ x_2 \in \mathbb{T}^1, \br
	\vu_{0,\delta} (x_1, x_2) &= \vu_\ell \ \mbox{if}\ x^i_1 - 2 \delta \leq x_1 \leq x^i_1 - \delta ,\ 
\vu_{0,\delta} (x_1, x_2) = \vu_r \ \mbox{if}\ x^i_1 + \delta \leq x_1 \leq x^i_1 + 2 \delta,\ x_2 \in \mathbb{T}^1
\label{P6}
	\end{align}
for $i=1, \dots,N$.

Next, we consider $\delta > 0$ small so that $(\vr_{0,\delta}, \vt_{0, \delta}, \vu_{0, \delta})$ defined in \eqref{P4}, \eqref{P6} may be extended to $\TD$ as 
$(\widetilde{\vr}_{0, \delta}, \widetilde{\vt}_{0, \delta}, \widetilde{\vu}_{0, \delta})$, 
\begin{equation} \label{P3}
	\widetilde{\vr}_{0, \delta} \in W^{k,2}(\TD),\ \inf_{\TD} \widetilde{\vt}_{0, \delta} > 0,\ \widetilde{\vt}_{0,\delta} \in W^{k,2}(\TD),\ \inf_{\TD} \widetilde{\vt}_{0,\delta} > 0,\ \widetilde{\vu}_{0, \delta} \in W^{k,2}(\TD; R^2),\ k > 2.	
\end{equation}

Finally, we introduce the initial data $(\vr_{0,\ep}, \vt_{0, \ep}, \vu_{0, \ep})$ supplementing
\eqref{P4}--\eqref{P6} by  
\begin{align} 
	\vr_{0,\delta} (x_1, x_2) &= \vr_\ell \ \mbox{if}\ x^i_1 - \delta \leq x_1 < x^i_1  ,\ 
	\vr_{0,\delta} (x_1, x_2) = \vr_r \ \mbox{if}\ x^i_1  \leq x_1 \leq x^i_1 + \delta,\ x_2 \in \mathbb{T}^1, \br
	\vt_{0,\delta} (x_1, x_2) &= \vt_\ell \ \mbox{if}\ x^i_1 - \delta \leq x_1 < x^i_1  ,\ 
	\vt_{0,\delta} (x_1, x_2) = \vt_r \ \mbox{if}\ x^i_1 \leq x_1 \leq x^i_1 + \delta,\ x_2 \in \mathbb{T}^1, \br
	\vu_{0,\delta} (x_1, x_2) &= \vu_\ell \ \mbox{if}\ x^i_1 -  \delta \leq x_1 < x^i_1 ,\ 
	\vu_{0,\delta} (x_1, x_2) = \vu_r \ \mbox{if}\ x^i_1 \leq x_1 \leq x^i_1 + \delta,\ x_2 \in \mathbb{T}^1
	\label{P7}
\end{align}
for $i=1, \dots,N$, and fixing $\delta = \delta(\ep)$ so that
\begin{equation} \label{P8} 
	\left\| \Big(\vr_{0,\delta(\ep)} - \vr_0;  \vt_{0,\delta(\ep)} - \vt_0; \vu_{0,\delta(\ep)} - \vu_0 \Big) \right\|_{L^q(\TD; R^4)} \leq \ep.
\end{equation}

From this moment on, the parameter $\delta = \delta(\ep)$ is fixed. Obviously, the data $(\vr_{0,\ep}, \vt_{0, \ep}, \vu_{0, \ep})=(\vr_{0,\delta(\ep)},  \vt_{0,\delta(\ep)}, \vu_{0,\delta(\ep)}) $ are piecewise 
smooth and satisfy \eqref{M3} in accordance with the conclusion of Theorem \ref{TM1}.

\subsection{Pasting solutions} \label{Ppasting}

Given the smooth data $(\widetilde{\vr}_{0, \delta}, \widetilde{\vt}_{0, \delta}, \widetilde{u}_{0, \delta})$ as in \eqref{P3}, the Euler system admits a (unique) local in time smooth 
solution $(\tvr, \tvt, \tvu)$ defined on a maximal time interval $T_{\rm max} > 0$. As regular solutions of the Euler system exhibit finite speed of propagation and the initial data 
satisfy \eqref{P6}, we deduce
\begin{align} 
	\tvr (t, x_1, x_2) &= \vr_\ell \ \mbox{if}\ x^i_1 - \frac{7}{4} \delta \leq x_1 \leq x^i_1 - \frac{5}{4} \delta ,\ 
	\tvr (t, x_1, x_2) = \vr_r \ \mbox{if}\ x^i_1 + \frac{5}{4} \delta \leq x_1 \leq x^i_1 + \frac{7}{4} \delta,\ x_2 \in \mathbb{T}^1, \br
	\tvt (t,x_1, x_2) &= \vt_\ell \ \mbox{if}\  x^i_1 - \frac{7}{4} \delta \leq x_1 \leq x^i_1 - \frac{5}{4} \delta ,\ 
	\tvt (t,x_1, x_2) = \vt_r \ \mbox{if}\ x^i_1 + \frac{5}{4} \delta \leq x_1 \leq x^i_1 + \frac{7}{4} \delta,\ x_2 \in \mathbb{T}^1, \br
	\tvu (t,x_1, x_2) &= \vu_\ell \ \mbox{if}\  x^i_1 - \frac{7}{4} \delta \leq x_1 \leq x^i_1 - \frac{5}{4} \delta  ,\ 
	\tvu(t,x_1, x_2) = \vu_r \ \mbox{if}\ x^i_1 + \frac{5}{4} \delta \leq x_1 \leq x^i_1 + \frac{7}{4} \delta ,\ x_2 \in \mathbb{T}^1
	\label{P9}
\end{align}
for $i=1, \dots,N$ for all $0< t \leq T_S$ whenever $T_S = T_S(\delta) > 0$ is small enough.

Next, consider the family of solutions $(\vr^n_{i,R}, \vt^n_{i,R}, \vu^n_{i,R})_{n=1}^\infty$ defined as 
\begin{align}
(\vr^n_{i,R}, \vt^n_{i,R}, \vu^n_{i,R})(t,x_1,x_2) = 
(\vr^n_{R}, \vt^n_{R}, \vu^n_{R})(t, x_1 - x^i_1, x_2),\ i=1, \dots, N,
\nonumber 	
	\end{align}
where $(\vr^n_{R}, \vt^n_{R}, \vu^n_{R})_{n=1}^\infty$ is the infinite family of distinct admissible solutions emanating from the Riemann data $
(\vr_\ell, \vt_\ell, \vu_\ell)$, $(\vr_r, \vt_r, \vu_r)$, the existence of which is guaranteed by Proposition \ref{PR1}. As they satisfy \eqref{R1a} with the 
same constant $\lambda$, there exists a positive time $T_R = T_R(\delta)$ such that 
\begin{align} 
(\vr^n_{R}, \vt^n_{R}, \vu^n_{R})(t, x_1, x_2) &= (\vr_\ell, \vt_\ell, \vu_\ell) \ \mbox{for} \ x_1 < -\delta,\ x_2 \in \mathbb{T}^1\br 
(\vr^n_{R}, \vt^n_{R}, \vu^n_{R})(t, x_1, x_2) &= (\vr_r, \vt_r, \vu_r)\ \mbox{for}\ x_1 > \delta,\ x_2 \in \mathbb{T}^1
\label{P10} 
\end{align}
for all $0 < t \leq T_R$. 

Finally, we set $T = \min \{ T_S, T_R \}$ and define a family of solutions emanating from the initial data $(\vr_{0,\ep}, \vt_{0, \ep}, \vu_{0, \ep})$ as follows
\begin{align} 
	(\vr^n, \vt^n, \vu^n)(t,x_1,x_2) &= (\tvr, \tvt, \tvu)(t,x_1,x_2) \ \mbox{if}\ 
|x_1 - x^i_1| \geq \frac{7}{4} \delta \ \mbox{ for all}\ i = 1, \dots, N \br	
(\vr^n, \vt^n, \vu^n)(t,x_1,x_2) &= (\vr^n_{i,R}, \vt^n_{i,R}, \vu^n_{i,R})(t,x_1,x_2) \ \mbox{if}\ |x_1 - x^i_1| \leq \frac{7}{4} \delta
\ \mbox{ for some}\ i = 1, \dots, N,
\label{P11}	
\end{align}
$x_2 \in \mathbb{T}^1$, $0 < t \leq T$. In view of \eqref{P9}, \eqref{P10}, the solutions 
$(\tvr, \tvt, \tvu)$ and $(\vr^n_R, \vt^n_R, \vu^n_R)$ coincide whenever 
\[
\frac{5}{4}\delta \leq |x_1 - x^i_1| \leq \frac{7}{4} \delta \ \mbox{for some}\ i = 1, \dots,N
\]
so \eqref{P11} yields a collection of admissible solutions satisfying \eqref{M2}, \eqref{M4}.

\section{Concluding remarks}
\label{D}

We have stated the main result in the spatial dimension $d=2$. As noticed by the authors of \cite{KlKrMaMa}, 
the ill posedness of the Riemann problems can be easily extended to $d=3$; whence the same applies to Theorem \ref{TM1}.

\subsection{Entropy producing solutions}
\label{EPS}

The method of proof of Theorem \ref{TM1} presented in Section \ref{P} can be slightly modified so to construct, starting from the initial data $(\vr_{0, \ep}, \vt_{0, \ep}, \vu_{0, \ep})$, infinitely many admissible solutions $(\vr^n, \vt^n, \vu^n)_{n \in N}$ which are moreover producing entropy.
Indeed, we can proceed as follows:
\begin{itemize}
\item we first fix arbitrarily an $i^\ast$ with $1\leq i^\ast \leq N$ and correspondingly select the point $x_1^{i^\ast}$ among the $N$ distinct points in \eqref{P1};
\item we proceed with the approximation of the initial data exactly as in Section \ref{Papprox};
\item we paste the solutions as performed in Section \ref{Ppasting} for all $1\leq i \leq N$, $i \neq i^\ast$;
\item for $i=i^\ast$, we define 
\begin{align}
(\vr_{i^\ast,R}, \vt_{i^\ast,R}, \vu_{i^\ast,R})(t,x_1,x_2) = 
(\vr_{RS}, \vt_{RS}, \vu_{RS})(t, x_1 - x^{i^\ast}_1, x_2),
\nonumber 	
	\end{align}
where $(\vr_{RS}, \vt_{RS}, \vu_{RS})$ is the standard self-similar Riemann solution emanating from the 
Riemann data $(\vr_\ell, \vt_\ell, \vu_\ell)$, $(\vr_r, \vt_r, \vu_r)$ and we observe that such a self-similar solution will contain a shock due to the choice of the Riemann initial data in Proposition \ref{PR1} (see \cite[Theorem 1.1]{KlKrMaMa});
we also notice that $(\vr_{i^\ast,R}, \vt_{i^\ast,R}, \vu_{i^\ast,R})(t,x_1,x_2)$ satisfy \eqref{P10};
\item
finally, we define the family of solutions
\begin{align} 
	(\vr^n, \vt^n, \vu^n)(t,x_1,x_2) &= (\tvr, \tvt, \tvu)(t,x_1,x_2) \ \mbox{if}\ 
|x_1 - x^i_1| \geq \frac{7}{4} \delta \ \mbox{ for all}\ i = 1, \dots, N \br	
(\vr^n, \vt^n, \vu^n)(t,x_1,x_2) &= (\vr^n_{i,R}, \vt^n_{i,R}, \vu^n_{i,R})(t,x_1,x_2) \ \mbox{if}\ |x_1 - x^i_1| \leq \frac{7}{4} \delta
\br &\mbox{ for some}\ i = 1, \dots, N, \ i\ne i^\ast, \br
(\vr^n, \vt^n, \vu^n)(t,x_1,x_2) &= (\vr_{i^\ast,R}, \vt_{i^\ast,R}, \vu_{i^\ast,R})(t,x_1,x_2) \ \mbox{if}\ |x_1 - x^{i^\ast}_1| \leq \frac{7}{4} \delta
\nonumber
\end{align}
$x_2 \in \mathbb{T}^1$, $0 < t \leq T$, with $T$ choosen as in Section \ref{P}, thus obtaining an infinite family of entropy producing admissible solutions emanating from the initial data $(\vr_{0, \ep}, \vt_{0, \ep}, \vu_{0, \ep})$.
\end{itemize}

\subsection{Problems on general domains}

Our method can be used to show the density of wild data on a general domain $\Omega \subset \mathbb{R}^2$, where the Euler system is supplemented with the impermeability boundary condition 
\begin{equation} \label{P12}
	\vu\cdot \vc{n}|_{\partial \Omega} = 0.
	\end{equation}
Indeed the construction is easy to adapt as soon as the domain contains a strip $(a,b) \times (c,d)$, specifically 
\begin{equation} \label{P13}
(a,b) \times (c,d) \subset \Omega, \ [a,b] \times \{ c \},\ [a,b] \times \{ d \} \subset \partial \Omega.
\end{equation}	

\begin{itemize}

\item	
The necessary local existence result was proved by Schochet \cite{SCHO1}.

\item To ensure the satisfaction of the compatibility conditions, we may consider the initial data \eqref{M1}, with $\vr_0$, $\vt_0$ constant and $\vu_0 = 0$ in a small neighbourhood of $\partial \Omega$.  

\item The construction of the modified data on a single strip $(a,b) \times (c,d)$ can be carried out exactly as in Section \ref{Papprox} yielding the data with vanishing second component $u_{0,2}$ in the 
strip $(a,b) \times (c,d)$, specifically satisfying the impermeability condition \eqref{P12} as well as the necessary compatibility conditions.

\end{itemize}

The above delineated construction yields the same conclusion as Theorem \ref{TM1} except the property \eqref{M2}. In particular, Corollary \ref{CM1} remains valid. The entropy producing solutions 
claimed in Section \ref{EPS} can be obtained in the same way.

\def\cprime{$'$} \def\ocirc#1{\ifmmode\setbox0=\hbox{$#1$}\dimen0=\ht0
	\advance\dimen0 by1pt\rlap{\hbox to\wd0{\hss\raise\dimen0
			\hbox{\hskip.2em$\scriptscriptstyle\circ$}\hss}}#1\else {\accent"17 #1}\fi}


\begin{thebibliography}{10}
	
	\bibitem{ABKlKrMaMa}
	H.~Al~Baba, C.~Klingenberg, O.~Kreml, V.~M\'{a}cha, and S.~Markfelder.
	\newblock Nonuniqueness of admissible weak solution to the {R}iemann problem
	for the full {E}uler system in two dimensions.
	\newblock {\em SIAM J. Math. Anal.}, {\bf 52}(2):1729--1760, 2020.
	
	\bibitem{BenSer}
	S.~Benzoni-Gavage and D.~Serre.
	\newblock {\em Multidimensional hyperbolic partial differential equations,
		{F}irst order systems and applications}.
	\newblock Oxford Mathematical Monographs. The Clarendon Press Oxford University
	Press, Oxford, 2007.
	
	\bibitem{BuDeSzVi}
	T.~Buckmaster, C.~de~Lellis, L.~Sz\'{e}kelyhidi, Jr., and V.~Vicol.
	\newblock Onsager's conjecture for admissible weak solutions.
	\newblock {\em Comm. Pure Appl. Math.}, {\bf 72}(2):229--274, 2019.
	
	\bibitem{BucVic}
	T.~Buckmaster and V.~Vicol.
	\newblock Convex integration and phenomenologies in turbulence.
	\newblock {\em EMS Surv. Math. Sci.}, {\bf 6}(1):173--263, 2019.
	
	\bibitem{ChVaYu}
	R.~M. Chen, A.~F. Vasseur, and Ch. Yu.
	\newblock Global ill-posedness for a dense set of initial data to the
	isentropic system of gas dynamics.
	\newblock {\em Adv. Math.}, 393:Paper No. 108057, 46, 2021.
	
	\bibitem{ChiDelKre}
	E.~Chiodaroli, C.~{D}e {L}ellis, and O.~Kreml.
	\newblock Global ill-posedness of the isentropic system of gas dynamics.
	\newblock {\em Comm. Pure Appl. Math.}, {\bf 68}(7):1157--1190, 2015.
	
	\bibitem{ChFe2022}
	E.~Chiodaroli and E.~Feireisl.
	\newblock On the density of {"}wild{"} initial data for the barotropic {E}uler
	system.
	\newblock {\em {\bf arxiv preprint No. 2208.04810}}, 2022.
	
	\bibitem{ChiKre}
	E.~Chiodaroli and O.~Kreml.
	\newblock On the energy dissipation rate of solutions to the compressible
	isentropic {E}uler system.
	\newblock {\em Arch. Ration. Mech. Anal.}, {\bf 214}(3):1019--1049, 2014.
	
	\bibitem{D4a}
	C.~M. Dafermos.
	\newblock {\em Hyperbolic conservation laws in continuum physics}, volume 325
	of {\em Grundlehren der Mathematischen Wissenschaften [Fundamental Principles
		of Mathematical Sciences]}.
	\newblock Springer-Verlag, Berlin, fourth edition, 2016.
	
	\bibitem{FeKlKrMa}
	E.~Feireisl, C.~Klingenberg, O.~Kreml, and S.~Markfelder.
	\newblock On oscillatory solutions to the complete {E}uler system.
	\newblock {\em J. Differential Equations}, {\bf 269}(2):1521--1543, 2020.
	
	\bibitem{FeKlMark}
	E.~Feireisl, C.~Klingenberg, and S.~Markfelder.
	\newblock On the density of ``wild'' initial data for the compressible {E}uler
	system.
	\newblock {\em Calc. Var. Partial Differential Equations}, {\bf 59}(5):Paper
	No. 152, 17, 2020.
	
	\bibitem{KlKrMaMa}
	C.~Klingenberg, O.~Kreml, V.~M\'{a}cha, and S.~Markfelder.
	\newblock Shocks make the {R}iemann problem for the full {E}uler system in
	multiple space dimensions ill-posed.
	\newblock {\em Nonlinearity}, {\bf 33}(12):6517--6540, 2020.
	
	\bibitem{MarKli}
	S.~Markfelder and C.~Klingenberg.
	\newblock The {R}iemann problem for the multidimensional isentropic system of
	gas dynamics is ill-posed if it contains a shock.
	\newblock {\em Arch. Ration. Mech. Anal.}, {\bf 227}(3):967--994, 2018.
	
	\bibitem{SCHO1}
	S.~Schochet.
	\newblock The compressible {E}uler equations in a bounded domain: {E}xistence
	of solutions and the incompressible limit.
	\newblock {\em Commun. Math. Phys.}, {\bf 104}:49--75, 1986.
	
	\bibitem{SzeWie}
	L.~Sz{\'e}kelyhidi and E.~Wiedemann.
	\newblock Young measures generated by ideal incompressible fluid flows.
	\newblock {\em Arch. Rational Mech. Anal.}, {\bf 206}:333--366, 2012.
	
\end{thebibliography}

\end{document}